\documentclass[12pt]{article}
\usepackage[utf8]{inputenc}

\usepackage[top=1in, bottom=1in, left=1.25in, right=1.25in]{geometry}
\usepackage{graphicx} 
\usepackage{subfig}
\usepackage{amsthm}
\usepackage[comma,numbers]{natbib}
\usepackage{ulem}
\usepackage{color}
\usepackage{lineno}  
\usepackage{hyperref}
\usepackage{amsmath}
\usepackage{amssymb}
\usepackage{indentfirst} 
\usepackage{mathrsfs}

\numberwithin{equation}{section}

\newcommand{\eps}{\varepsilon}

\newcommand\keywords[1]{\textbf{Keywords}: #1}

\theoremstyle{plain}

\numberwithin{equation}{section}

\title{An approximate solution of a perturbed Fokker-Planck equation}
\author{\small{Yan Luo$^{2}$, Kaicheng Sheng$^1$\footnote{Corresponding author. E-mail addresses:  {\it k.sheng@sdu.edu.cn} (K. Sheng)}} \\
\small{$^1$ School of Mathematics and Statistics, }\\ \small{Shandong University, Weihai, 264209, China }
\\
\small{$^2$ Research Centre for Mathematics and Interdisciplinary Sciences,}\\ \small{Shandong University, Qingdao, 266237, China }}

\date{}

\begin{document}

\maketitle
 

\begin{abstract} 

This paper focuses on finding an approximate solution of a kind of Fokker-Planck equation with time-dependent perturbations. A formulation of the approximate solution of the equation is constructed, and then the existence of the formulation is proved. The related Hamiltonian dynamical system explains the estimations. Examples of the Ornstein-Uhlenbeck process model and the nonlinear Langevin equation are used to validate the proposed results. Our work provides a more comprehensive understanding of the long-time behaviour of systems described by this Fokker-Planck equation and the corresponding stochastic differential equation.

\bigskip

\noindent\keywords{Fokker-Planck equation, perturbation, approximate solution, existence of the solution.}
\end{abstract}

\section{Introduction}

Stochastic Differential Equations (SDEs) are powerful tools for modelling systems with inherent randomness and noise. SDEs are widely used to formulate stochastic systems \cite{klebaner2005, oksendal2003}.  While SDEs describe the system's evolution, it is challenging to analyse the corresponding probability distribution directly. Probability density functions play an important role in calculating the probability. The Fokker-Planck equation (FPE) offers a complementary approach by focusing on the probability density function  \cite{risken1989}. 

The FPE arises in many fields, such as physics and chemistry \cite{vankampen2007}, biology \cite{mcadams1997}, finance \cite{black1973}, and engineering \citep{khasminskii2012}. Many studies have employed numerical methods and approximation techniques to address different kinds of FPEs \cite{fox1986, higham2001, kloeden1992, risken1996, xu2019}. While numerical studies are widely used to obtain solutions of the FPEs, analytical solutions of some FPEs provide more valuable insights into dynamics and statistical properties. The existence of the solutions of many FPEs has been proved \cite{bogachev2015}. Analytical methods for solving certain simple FPEs are discussed in a series of books, such as \cite{chung2001, risken1984}, but analytical solutions of more general FPEs are difficult to obtain.

Significant attention has been devoted to studying the behaviour of systems described by the Fokker-Planck equation with time-dependent perturbations. These perturbations can arise from various sources, including external fields, interactions with other particles, or changes in the system's internal dynamics. Understanding the perturbations' impact on the probability distribution's evolution is important \cite{gardiner2004, pavliotis2014}. Variational methods have been developed to find perturbative solutions to the Fokker-Planck equation with nonlinear drift \cite{Dreger2005}, and an approximation of the solution in the Fokker-Planck equation for Brownian motion with time-varying control parameters is given \cite{Wadia2022}. Perturbation technique was developed to approximate the corresponding Fokker–Planck operator \cite{Frank2016, Heninger2018}. Perturbation series are used to approach the solution of the Fokker-Planck equation \cite{Barcilon1996, Ma1992}. Analytical solutions for specific cases of time-dependent perturbed FPEs are studied \cite{ caldirola1951, kanai1950}. Some numerical approximations for the perturbed Fokker-Planck equation are given \cite{Narayanan2012, Odibat2007}.

An important example of studying the SDEs is the Benzi-Parisi-Sutera-Vulpiani (BPSV) stochastic resonance model of the dynamics of Quaternary glaciation, which is a theoretical framework that describes a phenomenon where the presence of internal or external noise in a nonlinear system can enhance the system's output response \cite{gritsyk1995, kleinen1996, sutera1983}. One important way to understand its probability distribution is to study the probability density function from FPE. Most people solve the FPEs by numerical methods. We generalised the BPSV model to a kind of stochastic differential equation and gave an approximate solution of the corresponding FPE.

In this paper, we focus on finding an approximate solution of a representative FPE with time-dependent perturbations. We construct a formulation of the approximate solution of the equation and then prove the existence of the formulation. We prove that the solution of the corresponding non-perturbed equation can approximately evaluate the solution of the studied FPE in the long-time case. An explanation of the estimations is detailed in the related Hamiltonian dynamical system. Examples of the Ornstein-Uhlenbeck process model and the nonlinear Langevin equation are taken to verify our results. Our work will provide a more comprehensive understanding of the long-time behaviour of systems described by this FPE.

\section{Deriving the Fokker-Planck equation}

This paper investigates the long-time behaviour of a 1-dimensional stochastic system governed by a stochastic differential equation (SDE) of the form
\begin{equation}
\label{SDE}
    dX_t=b\left(t,X_t\right)dt+\sigma dW_t,
\end{equation}
where $X_t\in \mathbb{R}$, $W_t$ is a standard Wiener process, $b\left(t,x\right): \mathbb{R}^{+} \times \mathbb{R} \to  \mathbb{R}$ is the drift coefficient function, $\sigma\in \mathbb{R}$ is the positively defined diffusion constant characterizing the strength of the stochastic noise and it is assumed to be bounded. 

The system is subject to both deterministic and random forces, with the drift term encapsulating the influence of deterministic dynamics and the diffusion term representing the stochasticity of the system. To understand the evolution of the system's distribution, it is important to find the probability density function (PDE) where its integral for an interval provides the probability of a value occurring in that interval.  We employ the Fokker-Planck equation (FPE) to provide a deterministic description of the evolution of the PDF $p(t, x)$. The corresponding FPE for this equation (\ref{SDE}) is given by:
\begin{equation}
\label{FPE}
\frac{\partial p(t, x)}{\partial t} = -\frac{\partial}{\partial x} \left[ b(t, x) p(t, x) \right] + \frac{1}{2} \sigma^{2} \frac{\partial^2 p(t, x)}{\partial x^2} , 
\end{equation}    
where $p(t, x)\in\mathbb{R}^+$ is the PDF of the system at time $t \in \mathbb{R}^{+}$ and position $x\in \mathbb{R}$, and the equation accounts for both the drift and diffusion terms governing the evolution of the density. 

The existence of solutions in this FPE has been rigorously established in previous works \cite{bogachev2011, bogachev2015} under appropriate conditions. Moreover, the existence of the periodic density function $p\left(t,x\right)$ in equation (\ref{FPE}) is given in \cite{chen2017, feng2023}. 

In this paper, $b\left(t,x\right)$ is supposed to be a separable function, i.e.
\begin{equation}
\label{sepera}
b\left(t,x\right)=-V'\left(x\right)+\eps h\left(t\right), \varepsilon>0,
\end{equation} 
where $V\left(x\right) \in \mathbb{R}$ is a potential function that depends only on the state variable $x$, $h\left(t\right) \in \mathbb{R}$ is a bounded time-dependent perturbation. The perturbation parameter $\eps$ is assumed to be small, which allows for the investigation of the system's behaviour in the limit of small perturbations. The system can be studied in more detail through an asymptotic approach that considers the long-time behaviour of the solution. The long-time behaviour of the system is considered by setting time $t \in \mathbb{T}=[\delta_s,\delta_l]$, $\delta_s>0$ and $\delta_l$ is large enough. 

In this framework, the FPE can be rewritten as:
\begin{equation}
\label{FPEB}
\frac{\partial p(t, x)}{\partial t} = -\frac{\partial}{\partial x} \left[ \left(-V'\left(x\right)+\eps h\left(t\right)\right) p(t, x) \right] + \frac{1}{2} \sigma^{2} \frac{\partial^2 p(t, x)}{\partial x^2}.
\end{equation}    

This paper aims to find an approximate solution of this equation for small values of $\eps$. Specifically, we will seek an approximate solution:
\begin{equation*}
\label{appsol}
\tilde{p}(t, x)=\hat{p}\left(x\right)+\eps\hat{p}\left(x\right)\kappa\left(t, x\right),
\end{equation*}  
where $\hat{p}(x)$ is the leading-order solution, and $\kappa(t, x)$ represents the first-order correction term. The solution $\hat{p}(x)$ corresponds to the solution of the FPE in the absence of the time-dependent perturbation ($\eps=0$), which is given by the equation:
\begin{equation}
\label{FPE2}
0= \frac{\partial}{\partial x} \left[V'(x)\,\hat{p}(x) \right] + \frac{1}{2} \sigma^{2} \frac{\partial^2 \hat{p}(x)}{\partial x^2}.
\end{equation}    
This is an autonomous Fokker-Planck equation whose solution
\begin{equation}
\label{hatp}
\hat{p}\left(x\right)=C \, e^{-\frac{2}{\sigma^{2}} V\left(x\right)}
\end{equation} 
is well-known and represents the stationary distribution of the system in the absence of time-dependent forcing \cite{bogachev2015}. For convenience, it is reasonable to take $C=1$ in the following. The full approximate solution $\tilde{p}(t, x)$, including the first-order correction, can then be expressed as:
\begin{equation}
\label{appsol2}
\tilde{p}(t, x)=e^{-\frac{2}{\sigma^{2}} V\left(x\right)} +\eps e^{-\frac{2}{\sigma^{2}} V\left(x\right)}\kappa\left(t, x\right),
\end{equation}  
which represents the probability density at time $t$ and position $x$ under the influence of the small time-dependent perturbation. We will prove that for large $t\in\mathbb{T}$, $\tilde{p}(t, x)$ satisfies
\begin{equation*}
\lim\limits_{\eps\to 0}\left| p(t, x)-\tilde{p}(t, x)\right | = 0, \quad  x\in\mathbb{R}.
\end{equation*} 
This approach provides a systematic method for obtaining the long-time asymptotic of the system, particularly in regimes where the perturbation is small.

\section{Finding function $\kappa(t, x)$}

\subsection{Existence of $\kappa(t, x)$}

To find the approximate solution (\ref{appsol2}), it is necessary to show the existence of the function $\kappa\left(t, x\right)$. In this section, details are given to show that $\kappa\left(t, x\right)$ is a solution of the following partial differential equation
\begin{equation}
\label{function_g}
\frac{\partial \kappa}{\partial t}=-V'\left(x\right) \frac{\partial \kappa}{\partial x}+\frac{1}{2} \sigma^{2} \frac{\partial^2 \kappa}{\partial x^2}+ \frac{2}{\sigma^2} V'\left(x\right) h\left(t\right) .
\end{equation}
We will show that $\kappa\left(t, x\right)$ is the desired function in the approximate solution (\ref{appsol2}) in the next section. 

To solve the PDE (\ref{function_g}), vector notations are used to express the equation more conveniently. We define the vector $\bar{\kappa}$ as
\begin{equation*}
\bar{\kappa}=\begin{pmatrix}
\kappa\\
\frac{\partial \kappa}{\partial x}\\
-V'\left(x\right)\kappa\\
-V''\left(x\right)\kappa
\end{pmatrix},
\end{equation*}
where $V'\left(x\right)$ and $V''\left(x\right)$ are assumed to be bounded. This allows us to work with the system that includes the function $\kappa(t, x)$ and its spatial derivatives. Next, we introduce the matrix operator $A$, which acts on the vector $\bar{\kappa}$ and incorporates all the differential terms. The operator $A$ is defined as
\begin{equation*}
A=\begin{pmatrix}
0 & -\frac{\sigma^2}{2}\frac{d}{dx} & -\frac{d}{dx} & 1\\
0 &-\frac{d}{dt}&0&0\\
0 & 0& -\frac{d}{dt}& 0\\
0 & 0& 0 & -\frac{d}{dt}
\end{pmatrix}.
\end{equation*}

Additionally, the perturbation due to the time-dependent function $h(t)$ is represented by the vector
\begin{equation*}
J=\begin{pmatrix}
h\left(t\right)\frac{2}{\sigma^2}V'\left(x\right)\\
0\\
0\\
0
\end{pmatrix}. 
\end{equation*}
This term reflects the driving force of the perturbation. With these definitions, the original equation (\ref{function_g}) for $\kappa(t, x)$ can be rewritten as
\begin{equation}
\label{function_g1}
-\frac{d\bar{\kappa}}{dt}=A\bar{\kappa}-J.
\end{equation}

Let $L$ be an operator  defined by 
$$
L\bar{\kappa}=-\frac{d\bar{\kappa}}{dt}.
$$
To solve this system, we use the resolvent operator. For $\hat{\lambda} \in \hat{\rho}(L)$, where $\hat{\rho}(L)$ is the resolvent set of the operator $L$, the resolvent $R(\hat{\lambda}, L)$ is given by
$$
R(\hat{\lambda}, L) = (\hat{\lambda} I - L)^{-1},
$$
and it acts on a function $v(t,x)$ as the following
$$
R(\hat{\lambda}, L) v(t,x) = \int_0^t e^{-\hat{\lambda}(t-s)} v(s,x) \, ds.
$$
This formula represents the evolution of $v(t,x)$ over time, with the exponential kernel $e^{-\hat{\lambda}(t-s)}$ capturing the time decay of the system. We usually take a large time $t$ to obtain a significant evolution. The resolvent operator satisfies the norm inequality
\begin{equation}
\label{RlambdaB}
\left\| R\left(\hat{\lambda},L\right)\right\|_{L^p}\leq \frac{1}{\hat{\lambda}},\quad \hat{\lambda}\in \mathbb{R},\, \hat{\lambda}>0.
\end{equation}
This inequality ensures that the resolvent operator is well-behaved as \( \hat{\lambda} \) increases, which is crucial for obtaining a solution of equation (\ref{function_g1}).

Then, we consider the Yosida approximation of the operator $L$, which is defined by
$$
L_n = n^2 R(n, L) - nI, \quad n\in \mathbb{N}.
$$
It is known from \cite{daprato1987} that as $n \to \infty$, the sequence converges to the operator $L$ in the $L^p$-norm:
\begin{equation}
\label{Ln}
    \lim_{n \to \infty} \| L_n v - L v \|_{L^p} = 0.
\end{equation}
Thus, the Yosida approximation provides a way to approximate the solution of the system as $n$ becomes large.

We now find a sequence \( v_n \in L^p \) that satisfies the equation
\begin{equation}
\label{bnvn}
-L_n\,v_n+A\,v_n-J=0, \quad n\in\mathbb{N}^{+}.
\end{equation}
In fact, this equation can be rewritten as
\begin{equation*}
\label{bnvn2}
\left(A-n\right)v_n=n^2\int_{0}^{t}e^{-n\left(t-s\right)}v_nds +J-2nv_n.
\end{equation*}
By taking the operator 
$$
(n I- A)^{-1} = R(n, A),
$$
then we obtain
\begin{equation}
\label{n2e}
   \left(-I+2nR\left(n,A\right)\right) v_n=n^2e^{-nt}\int_{0}^{t}e^{ns}R\left(n,A\right) v_n ds+R\left(n,A\right) J.
\end{equation}
Substituting
$$
w_n=\int_{0}^{t}e^{ns}R\left(n,A\right)v_n ds
$$
and 
$$
w_n'=e^{nt}R\left(n,A\right)v_n,
$$
into equation (\ref{n2e}), we have a new system
\begin{equation*}
\left(-I+2nR\left(n,A\right)\right)w_n'=n^2R\left(n,A\right)w_n+R^2\left(n,A\right)e^{nt} J,
\end{equation*}
as well as 
\begin{equation}
\label{wn}
\left(-I+2nR\left(n,A\right)\right)w_n=\int_{0}^{t}e^{n^2 R\left(n,A\right)\left(t-s\right)}R^2\left(n,A\right)e^{ns} J\,ds.
\end{equation}
This system can be solved iteratively, allowing us to obtain approximations for $v_n$, which will be proved to converge to the desired solution $\kappa(t, x)$ as $n \to \infty$.

To ensure that the sequence $v_n$ behaves well, we apply a bound on the norm of the resolvent operator $R(n, A)$. Let $\sigma$ be a finite positive constant. Then by the definition of $A$, we obtain
$$
nI-A=\begin{pmatrix}
nI & nI+\frac{d}{dx}\frac{\sigma^2}{2} & nI+\frac{d}{dx} & \left(n-1\right)I\\
0 & nI+\frac{d}{dt} &0 &0\\
0 & 0& nI+\frac{d}{ds} &0\\
0 &0 &0 & nI+\frac{d}{ds}
\end{pmatrix},
$$
and its reverse $R(n, A)=\left(nI-A\right)^{-1}$ is
\begin{equation}
\label{RNA}
\begin{aligned}
&\frac{\left(nI+\frac{d}{dt}\right)^{-1}}{n}\begin{pmatrix}
\left(nI+\frac{\sigma^2}{2}\frac{d}{dt}\right) & -\left(nI+\frac{\sigma^2}{2}\frac{d}{dx}\right) &-\left(nI+\frac{d}{dt}\right)& -\left(n-1\right)I\\
0 & nI & 0 &0\\
0& 0& nI &0\\
0& 0&0 & nI  
\end{pmatrix}.
\end{aligned}
\end{equation}
By the inequality (\ref{RlambdaB}), we get 
$$
\| R(n, A) \|_{L^p} \leq \frac{1}{n},
$$
as well as
$$
\left\|nR\left(n,A\right)\right\|_{L^p}<\infty, 
$$
and
$$
\left\|\left(-I+2nR\left(n,A\right)\right)^{-1}\right\|_{L^p}<\infty. 
$$
Then we have
\begin{equation*}
    \begin{aligned}
    &\left\|e^{n^2R\left(n,A\right)\left(t-s\right)}R^2\left(n,A\right)\right\|_{L^p}
    \\
    =&\left\|\sum_{k=0}^{\infty}\frac{1}{k!}\left(n^2R\left(n,A\right)\left(t-s\right)\right)^kR^2\left(n,A\right)\right\|_{L^p}\\
=&\left\|\sum_{k=0}^{\infty}\frac{1}{k!}\left(n^{2k}R^{k+2}\left(n,A\right)\left(t-s\right)^k\right)\right\|_{L^p}\\
    \le &M\sum_{k=0}^{\infty}\frac{n^{2k}\left(t-s\right)^k}{k!n^{k+2}}=\frac{Me^{n\left(t-s\right)}}{n^2},
    \end{aligned}
\end{equation*}
where
$$
M=\sup_{k\in \mathbb{N}, n>0}\left\|n^kR\left(n,A\right)^k\right\|_{L^p}<\infty.
$$ 
Thus 
\begin{equation*}
    \begin{aligned}
\left\|\left(-I+2nR\left(n,A\right)\right)v_n\right\|_{L^p} \le& M\left(\int_0^t\left\|J\right\|_{L^p} ds+\frac{1}{n}\left\|J\right\|_{L^p} \right)
\\
\le& MK\left\|J\right\|_{L^p},
    \end{aligned}
\end{equation*}
where $K$ is some finite positive constant. Specifically, we have
\begin{equation}
\label{vnlp}
    ||v_n||_{L^p}\le \bar{M}||J||_{L^p}
\end{equation} 
for some finite positive constant $\bar{M}$  depending on the properties of $R(n, A)$. As we assumed previously, $h$ and $V'$ are bounded, then $J$ is bounded. Moreover, inequality \ref{vnlp} ensures the boundedness of the sequence $v_n$ in the \( L^p \) -norm. This bound guarantees that the sequence does not grow to infinity, which is essential for the existence of a solution. Thus, we found a $L^{p}$-bounded sequence $v_n$ in equation (\ref{bnvn}).

Now we rewrite equation (\ref{bnvn}) as 
\begin{equation}
\label{Bnvn}
    -L_n\left(v_n-\bar{\kappa}\right)+A\left(v_n-\bar{\kappa}\right)-L_n\bar{\kappa}+L\bar{\kappa}=0,
\end{equation}
and define 
$$
\bar{\lambda}=L_n\bar{\kappa}-L\bar{\kappa},\quad V=v_n-\bar{\kappa}.
$$
By inequality (\ref{vnlp}) and equation (\ref{Bnvn}), we have
\begin{equation*}
\label{VLP}
    \begin{aligned}
    \left\|V\right\|_{L^p}=\left\|v_n-\bar{\kappa}\right\|_{L^p}\le \bar{M}\left\|\bar{\lambda}\right\|_{L^p}=\bar{M}\left\|L_n\bar{\kappa}-L\bar{\kappa}\right\|_{L^p},
    \end{aligned}
\end{equation*}
then by the formula (\ref{Ln}), we obtain 
\begin{equation}
\label{VLP2}
\lim\limits_{n\to \infty}\left\|v_n-\bar{\kappa}\right\|_{L^p}=0. 
\end{equation}
Therefore, we proved that $\kappa\left(t,x\right)$ exists as a solution of the PDE (\ref{function_g}) and it is bounded in the $L^p$-norm.

\subsection{Estimation of $\kappa(t, x)$}

The explicit solution of the partial differential equation (\ref{function_g}) is difficult to obtain analytically. Thus, an adequate approximation is considered in this subsection. From the convergence result in equation (\ref{VLP2}), we can express the function $\kappa\left(t,x\right)$ as the limit of the first component of $v_n$ in the $L^p$-norm.

It is necessary to calculate $v_n$. There are several methods to express $v_n$ explicitly. Directly, substituting equation (\ref{wn}) into equation (\ref{n2e}), we obtain the expression of $v_n$ as 
\begin{equation*}
    \begin{aligned}
&\left(-I+2nR\left(n,A\right)\right)v_n
\\
=&\left(-I+2nR\left(n,A\right)\right)^{-1}n^2e^{-nt}\int_{0}^{t}e^{n^2 R\left(n,A\right)\left(t-s\right)}R^2\left(n,A\right)e^{ns} J\,ds+R\left(n,A\right)J,
 \end{aligned}
\end{equation*}
where $R(n,A)$ is given by equation (\ref{RNA}). One may obtain the expression of $v_n$ with specific formulations of $V(x)$ and $h(t)$. This method sometimes brings a complicated calculation. In fact, recall that
$$
R(n, L) = (n I - L)^{-1},
$$
then 
$$
L_n^{-1} = (n^2 (n I - L)^{-1} - n I)^{-1}.
$$
By equation (\ref{bnvn}), 
$$
v_n = L_n^{-1}(A v_n - J)=(n^2 (n I - L)^{-1} - n I)^{-1} (A v_n - J).
$$ 
Therefore, $v_n$ is expressed explicitly:
$$
v_n = \left[(n^2 (n I - L)^{-1} - n I)^{-1} A - I \right]^{-1} J.
$$
Thus, we approximate $\kappa(t, x)$ by the sequence $v_n$ for sufficiently large $n$, with the convergence guaranteed in the $L^p$-norm. The approximation captures the essential features of the solution of the partial differential equation and can be used in further theoretical investigations.

\section{Construction of the approximate solution}

To estimate the exact solution \( p(t, x) \) of the FPE (\ref{FPEB}), we aim to derive an approximate solution \( \tilde{p}(t, x) \). We first recall the solution $\hat{p}(x)$, which satisfies the equation (\ref{FPE2}). In this case, \( \hat{p}(x) \) serves as the solution of the Fokker-Planck equation in the absence of time-dependent forces. 

We now seek an approximate solution \( \tilde{p}(t, x) \) that accounts for the time dependence introduced by the external forcing term $h(t)$. To do so, we begin by multiplying equation (\ref{function_g}) by $\hat{p}(x)$, yielding the equation
$\hat{p}\left(x\right)$ on both sides, 
\begin{equation}
\label{function_g3}
\hat{p}\left(x\right)\frac{\partial \kappa}{\partial t}=- V'\left(x\right)\hat{p}\left(x\right)\frac{\partial \kappa}{\partial x}+\frac{1}{2} \sigma^{2} \frac{\partial^2 \kappa}{\partial x^2}\hat{p}\left(x\right) + \frac{2}{\sigma^2} V'\left(x\right)\hat{p}\left(x\right) h\left(t\right) .
\end{equation}
where $\kappa(t, x)$ represents the correction term introduced by the perturbation. From equation (\ref{hatp}), we have the relation
\begin{equation}
\label{dhatq}
\frac{\partial \hat{p}\left(x\right)}{\partial x}=-\frac{2}{\sigma^2} V'\left(x\right)\hat{p}\left(x\right).
\end{equation}
Substituting this expression into the equation (\ref{function_g3}), we obtain
\begin{equation*}
\hat{p}\left(x\right)\frac{\partial \kappa}{\partial t}=-V'\left(x\right)\hat{p}\left(x\right)\frac{\partial \kappa}{\partial x}+\frac{1}{2} \sigma^{2} \frac{\partial^2 \kappa}{\partial x^2}\hat{p}\left(x\right) -\frac{\partial \hat{p}\left(x\right)}{\partial x} h\left(t\right), 
\end{equation*}
which can be rewritten as 
\begin{equation*}
\hat{p}\left(x\right)\frac{\partial \kappa}{\partial t}=V'\left(x\right)\hat{p}\left(x\right)\frac{\partial \kappa}{\partial x}+\frac{1}{2} \sigma^{2} \frac{\partial^2 \kappa}{\partial x^2}\hat{p}\left(x\right) -\frac{\partial \hat{p}\left(x\right)}{\partial x} h\left(t\right)-2V'\left(x\right)\hat{p}\left(x\right)\frac{\partial \kappa}{\partial x}. 
\end{equation*}
From equation (\ref{dhatq}), we have
\begin{equation}
    \label{FP2}
    \hat{p}\left(x\right)\frac{\partial \kappa}{\partial t}=V'\left(x\right)\hat{p}\left(x\right)\frac{\partial \kappa}{\partial x}+\frac{1}{2} \sigma^{2} \frac{\partial^2 \kappa}{\partial x^2}\hat{p}\left(x\right) -\frac{\partial \hat{p}\left(x\right)}{\partial x} h\left(t\right)+\sigma^2\frac{\partial \hat{p}}{\partial x}\frac{\partial \kappa}{\partial x}.
\end{equation}

Note that the approximate solution $\tilde{p}(t, x)=\hat{p}+\varepsilon\hat{p}k$ satisfies 
\begin{equation}
\label{tilde p}
    \frac{\partial \tilde{p}}{\partial t}=\frac{\partial \hat{p}}{\partial t}+\epsilon \frac{\partial \hat{p}}{\partial t}\kappa +\epsilon \hat{p} \frac{\partial \kappa}{\partial t},
\end{equation}
we put the equations (\ref{FPE2}) and (\ref{FP2}) into equation (\ref{tilde p}) to obtain an approximate solution $\tilde{p}(t, x)$ which satisfies the following:
\begin{equation}
\label{tildep}
\frac{\partial \tilde{p}(t, x)}{\partial t}=\frac{\partial}{\partial x} \left[ V'(x) \tilde{p}(t, x) \right]
+\frac{1}{2} \sigma^{2} \frac{\partial^2 \tilde{p}(t, x)}{\partial x^2}-\eps \frac{\partial \hat{p}\left(x\right)}{\partial x} h\left(t\right).
\end{equation}
This equation accounts for both the original dynamics of the Fokker-Planck equation and the perturbation introduced by the time dependence of $h(t)$. 

Let $q(t, x) = p(t, x) - \tilde{p}(t, x)$ represent the error between the exact solution $p\left(t,x\right)$ and the approximate solution $\tilde{p}(t, x)$. To derive the equation for $q(t, x)$, we combine equation (\ref{FPEB}) and equation (\ref{tildep}),  then
\begin{equation}
\label{FP3}
\frac{\partial q(t, x)}{\partial t}=\frac{\partial}{\partial x} \left[\left( V'(x) q(t, x) \right)\right]-\eps \left(\frac{\partial {p}\left(t,x\right)}{\partial x} -\frac{\partial \hat{p}\left(t, x\right)}{\partial x} \right) h\left(t\right)+\frac{1}{2} \sigma^{2} \frac{\partial^2 q(t, x)}{\partial x^2}.
\end{equation}

Meanwhile, from $\tilde{p}=\hat{p}+\epsilon \hat{p}\kappa$, we have 
 \begin{equation}
 \label{dgdx}
\eps\frac{\partial \tilde{p}\left(t, x\right)}{\partial x}=\eps \frac{\partial \hat{p}\left( x\right)}{\partial x}+\varepsilon^2 \frac{\partial \hat{p}\left(x\right)}{\partial x}\kappa\left(t, x\right)+\varepsilon^2 \frac{\partial \kappa\left(t, x\right)}{\partial x}\hat{p}\left(x\right).
 \end{equation}
Note that $p\left(t,x\right)=\tilde{p}\left(t,x\right)+q\left(t,x\right)$, then by equation (\ref{dgdx}), 
\begin{equation}
\begin{aligned}
    \label{phatp} \epsilon\left(\frac{\partial p}{\partial x}-\frac{\partial \hat{p}}{\partial x}\right)=&\epsilon\frac{\partial \tilde{p}}{\partial x}+\epsilon\frac{\partial q}{\partial x}\\
    =&\epsilon\frac{\partial q}{\partial x}+\varepsilon^2 \frac{\partial \hat{p}}{\partial x}\kappa+\varepsilon^2 \frac{\partial \kappa}{\partial x}\hat{p}.
    \end{aligned}
\end{equation}
Thus from equation (\ref{FP3})and equation (\ref{phatp}), we obtain the equation for $q(t, x)$:
\begin{equation}
\label{FP4}
\begin{aligned}
\frac{\partial q(t, x)}{\partial t}=& \frac{\partial}{\partial x} \left[ V(x) q(t, x) \right] +\frac{1}{2} \sigma^{2} \frac{\partial^2 q(t, x)}{\partial x^2}-\eps\frac{\partial q\left(t, x\right)}{\partial x}h\left(t\right)
\\
&-\eps^2 \frac{\partial \hat{p}\left(x\right)}{\partial x}\kappa\left(t, x\right)h\left(t\right)-\varepsilon^2 \frac{\partial \kappa\left(t, x\right)}{\partial x}\hat{p}\left(x\right) h\left(t\right).
\end{aligned}
\end{equation}

As we studied before, \( \bar{\kappa}(t, x) \) and $\hat{p}\left(x\right)$ are bounded for \( x \in  \mathbb{R} \). Specifically, 
$$
G\left(t,x\right)=\frac{\partial \hat{p}\left( x\right)}{\partial x}\kappa\left(t, x\right)h\left(t\right)+ \frac{\partial \kappa\left(t, x\right)}{\partial x}\hat{p}\left(x\right) h\left(t\right),
$$
is also bounded. Consequently, the equation (\ref{FP4}) for $q(t, x)$ can be written as
$$
-L\left(q\right)=\varepsilon^2\,G\left(t,x\right),
$$
where $L$ is an operator given by
\begin{equation*}
\begin{aligned}
L\left(q\right)=&\frac{\partial q(t, x)}{\partial t}-\frac{\partial}{\partial x} \left[ V(x) q(t, x) \right] -\frac{1}{2} \sigma^{2} \frac{\partial^2 q(t, x)}{\partial x^2}+\eps\frac{\partial q\left(t, x\right)}{\partial x}h\left(t\right).
\end{aligned}
\end{equation*}
Since $L\left(q\right)$ is bounded for any $x\in \mathbb{R},t\in \mathbb{T}$, then we have
$$
|p\left(t, x\right)-\tilde{p}\left(t, x\right)|\le \hat{C}\eps^2, \quad \, x\in\mathbb{R},\,t\in\mathbb{T}, 
$$
where $\hat{C}$ is some finite positive constant. Therefore, the error between $p(t, x)$ and $\tilde{p}(t, x)$ decays as $\varepsilon^2$, and we obtain the result
\begin{equation*}
\lim\limits_{\eps\to 0}\left| p(t, x)-\tilde{p}(t, x)\right | = 0, \quad  x\in\mathbb{R},\,t\in\mathbb{T}.
\end{equation*} 
Thus, the approximate solution $\tilde{p}(t, x)$ converges to the exact solution $p(t, x)$ as $\varepsilon \to 0$.

Moreover, the approximate solution $\tilde{p}(t, x)$ is given by
\begin{equation*}
\tilde{p}(t, x)=\hat{p}\left(x\right)+\eps\hat{p}\left(x\right) \kappa\left(t, x\right)=e^{-\frac{2}{\sigma^{2}} V\left(x\right)} +\eps e^{-\frac{2}{\sigma^{2}} V\left(x\right)}\kappa\left(t, x\right)
\end{equation*}  
which provides an approximate solution of the Fokker-Planck equation with accuracy $\eps$. 

For sufficiently small $\eps$ and large $t$, where $\left|V\left(x\right)\right|<c_1$ and $\left|\kappa\left(t, x\right)\right|< c_2$ (with $c_1$ and $c_2$ being constants), the error between $\tilde{p}(t, x)$ and $\hat{p}(x)$ is bounded by
\begin{equation*}
\begin{aligned}
\left|\tilde{p}(t, x)-\hat{p}(x)\right|=&\eps \left|e^{-\frac{2}{\sigma^{2}} V\left(x\right)}\kappa\left(t, x\right)\right| 
\\
< &\eps c_2  e^{\frac{2 c_1}{\sigma^{2}}} < c \eps, \quad  x\in\mathbb{R},\,t\in\mathbb{T},
\end{aligned}
\end{equation*}  
where $c=\max( c_2 e^{{2 c_1}/{\sigma^{2}}})$. Hence, the approximation error is of order $\eps$. Thus, the error between ${p}(t, x)$ from time-dependent FPE (\ref{FPEB}) and $\hat{p}(x)$ from the FPE in the absence of the time-dependent perturbation (equation (\ref{FPE2})) is bounded by
\begin{equation*}
\begin{aligned}
&\left| p(t, x)-\hat{p}(x)\right | 
\\
\leq &\left| p(t, x)-\tilde{p}(t, x)\right |+\left| \tilde{p}(t, x)-\hat{p}(x)\right | 
\\
= &\hat{C}\eps^2 + c \eps \leq \bar{c} \eps, \quad  x\in\mathbb{R},\,t\in\mathbb{T},
\end{aligned}
\end{equation*} 
where $\bar{c}$ is some finite positive constant.

In summary, for lower accuracy requirements, the solution ${p}(t, x)$ of the Fokker-Planck equation can be approximated by $\hat{p}(x)$ in the long-time behaviour. For higher accuracy, $\tilde{p}(t, x)$, the solution with the correction term, provides a better approximation. As $\eps \to 0$, ${p}(t, x)$ converges to $\hat{p}(x)$, and the approximation becomes increasingly accurate.

\section{Explanation in dynamical system}

In this section, we transform the Fokker-Planck equation (FPE) to its corresponding Hamilton-Jacobi form. To achieve this, we introduce an effective potential function $S(t, x)$ through the transformation:
\begin{equation*}
\label{FPEHJB}
p(t, x) = e^{-\frac{S(t,x)}{\sigma^2}},
\end{equation*} 
where $p(t, x)$ is the probability density function from FPE (\ref{FPE}). Then we obtain expressions for the first and second derivatives of $p(t, x)$ are 
\begin{equation*}
\frac{\partial p}{\partial t} = -\frac{1}{\sigma^2} e^{-\frac{S}{\sigma^2}} \frac{\partial S}{\partial t},
\end{equation*} 
\begin{equation*}
\frac{\partial p}{\partial x} = -\frac{1}{\sigma^2} e^{-\frac{S}{\sigma^2}} \frac{\partial S}{\partial x},
\end{equation*} 
and
\begin{equation*}
\frac{\partial^2 p}{\partial x^2} = -\frac{1}{\sigma^2} e^{-\frac{S}{\sigma^2}} \frac{\partial^2 S}{\partial x^2} + \frac{1}{\sigma^4} e^{-\frac{S}{\sigma^2}} \left( \frac{\partial S}{\partial x} \right)^2.
\end{equation*} 
Substituting these expressions into the FPE (\ref{FPE}) and factor out the non-zero term $e^{-{S}/{\sigma^2}}$. The equation is simplified to be
\begin{equation*}
\frac{\partial S}{\partial t} + \sigma^2 \frac{\partial b(t,x)}{\partial x} - b(t, x) \frac{\partial S}{\partial x} + \frac{\sigma^2}{2} \frac{\partial^2 S}{\partial x^2} - \frac{1}{2} \left( \frac{\partial S}{\partial x} \right)^2 = 0.
\end{equation*} 
This is the Hamilton-Jacobi equation for the system. We can then rewrite it in the standard form:
\begin{equation*}
\frac{\partial S}{\partial t} + H(t, x, \partial_x S, \partial_x^2 S) = 0,
\end{equation*} 
where the Hamiltonian $H$ is given by:
\begin{equation*}
H(t, x, \partial_x S, \partial_x^2 S) = -\sigma^2 \frac{\partial b(t,x)}{\partial x} + b(t,x) \frac{\partial S}{\partial x} - \frac{\sigma^2}{2} \frac{\partial^2 S}{\partial x^2} + \frac{1}{2} \left( \frac{\partial S}{\partial x} \right)^2.
\end{equation*} 
By substituting the expression for the potential $b(t,x)$ from equation (\ref{sepera}), we split the Hamiltonian into two parts:
\begin{equation*}
H(t, x, \partial_x S, \partial_x^2 S) = H_{0}(x, \partial_x S, \partial_x^2 S)+H_{pert}(t, \partial_x S),
\end{equation*}
where
\begin{equation*}
H_{0}(x, \partial_x S, \partial_x^2 S) = \sigma^2 V''(x) - V'(x) \frac{\partial S}{\partial x}- \frac{\sigma^2}{2} \frac{\partial^2 S}{\partial x^2} + \frac{1}{2} \left( \frac{\partial S}{\partial x} \right)^2
\end{equation*}
is the unperturbed Hamiltonian. The perturbed Hamiltonian is given by
\begin{equation*}
H_{pert}(t, \partial_x S) =  \varepsilon h(t) \frac{\partial S}{\partial x}.
\end{equation*}

Introducing a canonical momentum $y = \partial_x S$, and replacing $\partial_x^2 S$ as a derivative of $y$. We rewrite the Hamiltonian in terms of $y$ and $\partial_x y$: 
\begin{equation*}
H(t, x, y) = H_{0}(x, y)+H_{pert}(t, y),
\end{equation*}
where
$$
H_0(x, y) = \sigma^2 V''(x) - V'(x) y - \frac{\sigma^2}{2} \frac{\partial y}{\partial x} + \frac{1}{2} y^2.
$$
It is an autonomous function that depends on position $x$ and the canonical conjugate momentum $y$. The perturbed Hamiltonian become 
\begin{equation*}
H_{pert}(t, y) =  \varepsilon h(t) y.
\end{equation*}
The form of $V(x)$ determines whether periodic or quasi-periodic motions exist. Assuming the system exhibits periodic or quasi-periodic motion, according to \cite{arnold1974, akn}, the action-angle variable can be obtained by a canonical transformation $(x, y) \to (I, \theta)$. such that the action variable is computed as the integral over one period of the motion, the conjugate angle variable
$$
\theta = \frac{\partial S}{\partial I}.
$$
In the small perturbations $H_{pert}$,  the canonical momentum $y$ can be expanded as a Fourier series in $\theta$, thus the perturbed Hamiltonian becomes
\begin{equation*}
H_{pert}(I, \theta, t) = \varepsilon h(t) \sum_{k} C_k(I) e^{i k \theta},
\end{equation*}
where $C_k(I)$ are the Fourier coefficients of $y$ in terms of $\theta$. Thus the Hamiltonian 
\begin{equation*}
H(I, \theta, t) = H_0(I) + \varepsilon \sum_{k} h(t) C_k(I) e^{i k \theta}.
\end{equation*}

To analyse the persistence of periodic or quasi-periodic motion under the perturbation, we apply the Kolmogorov-Arnold-Moser (KAM) theory \citep{arnold1963proof, akn}. For an integrable Hamiltonian \( \mathcal{H}_0(\mathbf{I}) \), where \( \mathbf{I} = (I_1, \dots, I_n) \) are action variables, the frequencies are given by:  
$$
\omega(\mathbf{I}) = \nabla \mathcal{H}_0(\mathbf{I}) = \left( \frac{\partial \mathcal{H}_0}{\partial I_1}, \dots, \frac{\partial \mathcal{H}_0}{\partial I_n} \right).
$$
The theory gives two crucial conditions for the persistence of invariant tori under small perturbations:

\begin{itemize}
    \item[a)] Non-degeneracy condition: 
    
    The non-degeneracy condition requires that the frequency map \( \mathbf{I} \mapsto \omega(\mathbf{I}) \) is a local diffeomorphism, meaning the Hessian matrix is non-singular:  
    $$
    \det \left( \frac{\partial^2 \mathcal{H}_0}{\partial I_i \partial I_j} \right) \neq 0.
    $$ 
    \item[b)] Non-resonance condition (Diophantine condition):

    A frequency vector \( \omega(\mathbf{I}) \in \mathbb{R}^n \) satisfies the Diophantine condition if there exist constants \( C > 0 \) and \( \tau \geq n - 1 \), such that:  
    $$
    |\mathbf{k} \cdot \omega| \geq \frac{C}{\|\mathbf{k}\|^\tau} \quad \text{for all } \mathbf{k} \in \mathbb{Z}^n \setminus \{0\},
    $$
    where \( \|\mathbf{k}\| = |k_1| + \dots + |k_n| \).  
\end{itemize}

In a system with one degree of freedom ($n=1$), the non-degeneracy condition $\omega(I) = {\partial H_0} / {\partial I} \neq 0$ ensures that the frequencies associated with periodic or quasi-periodic motion are non-degenerate. The non-resonance condition is satisfied in the system naturally. The KAM theorem guarantees that for sufficiently small perturbations $H_{pert}(I, \theta, t)$, the majority of the periodic or quasi-periodic orbits of the unperturbed system $H_0(I)$  will persist. While the perturbation causes small shifts in the frequencies, the invariant tori corresponding to these orbits remain, though they may undergo slight deformations.

Finally, the transformation from the probability distribution \( p(t, x) \) to the function \( S(t, x) \) via the equation \( p(t, x) = e^{-{S(t, x)}/{\sigma^2}} \) is smooth and continuous. Therefore, in the stationary case, the solutions of the FPE (\ref{FPE}) approach those of the unperturbed, autonomous FPE (\ref{FPE2}). This provides a link between the stochastic properties of the system and its long-time dynamical behaviour.


\section{Applications of the equation}

\subsection{Ornstein-Uhlenbeck process model}

The Ornstein-Uhlenbeck process is an important extension of Brownian motion, considering the restoring force acting on particles (\cite{ornstein1930}). It is widely used to simulate the evolution of dynamical systems with random disturbances. The one-dimensional perturbed Ornstein-Uhlenbeck process is generated by the SDE
\begin{equation*}
\label{OUSDE}
    dX_t=\vartheta\left(\mu-X_t +\eps\cos t\right)dt+\sigma dW_t,
\end{equation*}
which can be regarded as an example of the SDE (\ref{SDE}) by taking 
$$
b\left(t, X_t\right)=\vartheta\left(\mu-X_t +\eps\cos t\right),
$$
where $\mu\in \mathbb{R}$ is the long-term mean or unconditional expectation, $\vartheta\in \mathbb{R}$ is the mean reversion rate, representing the speed at which the system reverts to the mean $\mu$. The corresponding Fokker-Planck equation is given by
\begin{equation}
\label{OUFPE}
\frac{\partial p(t, x)}{\partial t} = -\frac{\partial}{\partial x} \left[\vartheta\left(\mu-x +\eps\cos t\right) p(t, x) \right] + \frac{1}{2} \sigma^{2} \frac{\partial^2 p(t, x)}{\partial x^2} , 
\end{equation}
where $p(t, x)\in \mathbb{R}^{+}$ is the probability density function of the system at position $x\in \mathbb{R}$ and time $t \in \mathbb{T}$. The corresponding autonomous FPE is 
\begin{equation}
\label{FPE61}
0= -\frac{\partial}{\partial x} \left[\vartheta\left(\mu-x \right)\hat{p}(x) \right] + \frac{1}{2} \sigma^{2} \frac{\partial^2 \hat{p}(x)}{\partial x^2}
\end{equation}  
with the normalised solution 
\begin{equation*}
\begin{aligned}
\hat{p}(t, x)=&\frac{1}{\sqrt{\pi}\sigma}\exp\left[{\frac{2}{\sigma^{2}} \int_{0}^{x}\vartheta\left(\mu-s\right)ds}\right]  \\
=&\frac{1}{\sqrt{\pi}\sigma}\exp\left[{\frac{2}{\sigma^{2}} \vartheta\left(\mu x-\frac{1}{2}x^2\right)}\right].
\end{aligned}
\end{equation*}
Therefore, the approximate normalised solution of equation (\ref{OUFPE}) can be written as 
\begin{equation}
\label{OUapp}
\begin{aligned}
\tilde{p}(t, x)=&\frac{1}{\sqrt{\pi}\sigma}\exp\left[{\frac{2}{\sigma^{2}} \vartheta\left(\mu x-\frac{1}{2}x^2\right)}\right] 
\\
&+ \frac{\eps}{\sqrt{\pi}\sigma}\exp\left[{\frac{2}{\sigma^{2}} \vartheta\left(\mu x-\frac{1}{2}x^2\right)}\right]\kappa_1\left(t, x\right).
\end{aligned}
\end{equation}
where the function $\kappa_1\left(t, x\right)$ is bounded and can be estimated by equation 
\begin{equation*}
\label{function_g61}
\frac{\partial \kappa_{1}}{\partial t}=\vartheta\left(\mu-x \right) \frac{\partial \kappa_{1}}{\partial x}+\frac{1}{2} \sigma^{2} \frac{\partial^2 \kappa_{1}}{\partial x^2}- \frac{2}{\sigma^2} \vartheta\left(\mu-x \right) h\left(t\right) .
\end{equation*}

For $\vartheta=1$ and $\mu=0$, the solution (\ref{OUapp}) simplifies to 
\begin{equation}
\label{OUapp2}
\tilde{p}(t, x)=\frac{1}{\sqrt{\pi}\sigma}\exp\left(-\frac{x^2}{\sigma^{2}} \right) + \frac{\eps}{\sqrt{\pi}\sigma}\exp\left(-\frac{x^2}{\sigma^{2}} \right)\kappa_1\left(t, x\right), 
\end{equation}
and the explicit solution of equation (\ref{OUFPE}) is
\begin{equation}
\label{OUsol}
p(t, x)=\frac{1}{\sqrt{\pi\left(1-{\mathrm e}^{-2 t}\right)}\sigma}\cdot \exp\left(-\frac{\left( x+ \varepsilon Z_{1}\right)^{2}}{ \left(1-{\mathrm e}^{-2 t}\right) \sigma^{2}}\right),
\end{equation} 
where $Z_{1}=-\sin\left(t\right)-\cos\left(t\right)+{\mathrm e}^{-t}$ is bounded for $t\in \mathbb{T}$. 

When $t$ is large enough, $1-{\mathrm e}^{-2 t}$ can be approximately considered as $1$. Therefore, for large $t\in \mathbb{T}$, $x\in \mathbb{R}$, the error between the approximate solution (\ref{OUapp2}) and the explicit solution (\ref{OUsol}) is 
\begin{equation*}
\begin{aligned}
\delta(t, x)=&\left| p(t, x)-\tilde{p}(t, x)\right | 
\\
\approx& \left|\frac{1}{\sqrt{\pi}\sigma}\cdot \exp\left(-\frac{\left(2 x+ \varepsilon Z_{1}\right)^{2}}{4 \sigma^{2}}\right) -\frac{1}{\sqrt{\pi}\sigma}\exp\left(-\frac{x^2}{\sigma^{2}} \right) \right.
\\
& \left.- \frac{\eps}{\sqrt{\pi}\sigma}\exp\left(-\frac{x^2}{\sigma^{2}} \right)\kappa_1\left(t, x\right) \right|
\\
=&  \frac{1}{\sqrt{\pi}\sigma}\exp\left(-\frac{x^2}{\sigma^{2}} \right)\left|\exp\left(-\frac{4 \varepsilon x Z_{1}+ \varepsilon^2 Z_{1}^2}{4 \sigma^{2}}\right)-1-\eps \kappa_1\left(t, x\right)\right|   
\\
\approx& \frac{1}{\sqrt{\pi}\sigma}\exp\left(-\frac{x^2}{\sigma^{2}} \right)\left|-\frac{4 \varepsilon x Z_{1}+ \varepsilon^2 Z_{1}^2}{4 \sigma^{2}}-\eps \kappa_1\left(t, x\right)\right|   
\\
\sim& O(\eps).
\end{aligned}
\end{equation*} 

This reveals that our procedure gives a great approximation for some kinds of perturbed Fokker-Planck equations in long-time behaviour. 

\subsection{Nonlinear Langevin equation}

A nonlinear Langevin equation is a stochastic differential equation that describes the evolution of a physical system under the influence of both deterministic and random forces. The deterministic force in a nonlinear Langevin equation is a nonlinear function of the system's state. 

We consider the overdamped Langevin equation for a perturbed nonlinear oscillation system: 
\begin{equation*}
\label{LangSDE}
dY= \left[-a_1Y + a_2(\sin^2\frac{Y}{2} - a_3)\sin Y + \eps\cos t\right] dt + \sigma d W_t,
\end{equation*} 
where $Y \in \mathbb{R}$ is the angle variable in the system, $a_1$ represents the damping effect, which is the tendency of the system to resist motion or vibration, $a_2$ describes the nonlinear restoring force or nonlinear effects, and $a_3$ is a parameter within the nonlinear term, affecting the shape and strength of the nonlinear restoring force, $W_t$ is the Gaussian white noise. The force 
$$
F(t, y)=-a_1y + a_2(\sin^2\frac{y}{2} - a_3)\sin Y+ \eps\cos t
$$
refers to some nonlinear force with a deterministic time-dependent perturbation. The probability density function is determined by the following Fokker-Planck equation
\begin{equation}
\label{LangFPE}
\begin{aligned}
\frac{\partial p(t, y)}{\partial t} =& -\frac{\partial}{\partial y} \left[\left(-a_1y + a_2(\sin^2\frac{y}{2} - a_3)\sin y + \eps\cos t\right) p(t, y) \right] \\
&+ \frac{1}{2} \sigma^{2} \frac{\partial^2 p(t, y)}{\partial y^2}.
\end{aligned}
\end{equation} 
It is difficult to find the general analytical solution directly in this partial differential equation because of its nonlinearity, and it also admits a time-dependent perturbation, which is rather more complicated. However, if we take 
$$
F_0(y)=-a_1y + a_2(\sin^2\frac{y}{2} - a_3)\sin y, \quad F_1(t)=\cos t,
$$
the force $F(t, y)=F_0(y)+\eps F_1(t)$ satisfies the form of equation (\ref{sepera}). In this model, 
\begin{equation*}
\begin{aligned}
\hat{p}(t, x)=&\exp\left[{-\frac{2}{\sigma^{2}} \int_0^{y} F_0(\varphi)d\varphi}\right]  \\
=&\exp{ \left[ -\frac{2 a_{2}}{\sigma^{2}} \left(\sin\left(\frac{y}{2}\right)^{2}-a_{3}\right)^{2}+\frac{a_{1} y^{2}}{\sigma^{2}} \right]} 
\end{aligned}
\end{equation*}
is the solution of the corresponding autonomous FPE
\begin{equation*}
\label{FPE62}
0= -\frac{\partial}{\partial x} \left[\left(-a_1y + a_2(\sin^2\frac{y}{2} - a_3)\sin y \right)\hat{p}(x) \right] + \frac{1}{2} \sigma^{2} \frac{\partial^2 \hat{p}(x)}{\partial x^2}.
\end{equation*}  
Then the solution of the equation (\ref{LangFPE}) can be approximately given by
\begin{equation*}
\label{Langapp}
\begin{aligned}
\tilde{p}(t, y)=&\exp{ \left[ -\frac{2 a_{2}}{\sigma^{2}} \left(\sin\left(\frac{y}{2}\right)^{2}-a_{3}\right)^{2}+\frac{a_{1} y^{2}}{\sigma^{2}} \right]} 
\\
&+\eps \exp{ \left[ -\frac{2 a_{2}}{\sigma^{2}} \left(\sin\left(\frac{y}{2}\right)^{2}-a_{3}\right)^{2}+\frac{a_{1} y^{2}}{\sigma^{2}} \right]} \kappa_2\left(t, y \right),
\end{aligned}
\end{equation*}  
where the function $\kappa_2\left(t, x\right)$ is bounded and can be estimated by equation 
\begin{equation*}
\begin{aligned}
\frac{\partial \kappa_{2}}{\partial t}=&-\left(-a_1 y + a_2(\sin^2\frac{y}{2} - a_3)\sin y \right) \frac{\partial \kappa_{2}}{\partial y}
\\
&+\frac{1}{2} \sigma^{2} \frac{\partial^2 \kappa_{2}}{\partial y^2}- \frac{2}{\sigma^2} \left(-a_1y + a_2(\sin^2\frac{y}{2} - a_3)\sin y \right) h\left(t\right) .
\end{aligned}
\end{equation*}

\section{Conclusions}

In conclusion, this paper has focused on finding an approximate solution to a representative Fokker-Planck equation (FPE) with time-dependent perturbations. We constructed a formulation for the approximate solution and established its existence. Our analysis demonstrated that the solution of the corresponding non-perturbed equation can effectively approximate the solution of the studied FPE in the long-time behaviour. The estimations were further elucidated through the lens of the related Hamiltonian dynamical system. Examples of the Ornstein-Uhlenbeck process model and the nonlinear Langevin equation verified our results. This work contributes to a deeper understanding of the long-time behaviour of systems governed by such FPEs. It provides a foundation for future studies on perturbed systems and their asymptotic properties.

\bigskip

\section*{Acknowledgements}

The authors thank Prof. Huaizhong Zhao for the discussions. Yan Luo also thanks the National Natural Science Foundation of China (NSFC) for supporting this research (Grant: 12471142). 





\normalem
\bibliographystyle{plain}  
\bibliography{reference} 

\end{document}